\documentclass[a4paper,12pt,reqno]{amsart}

\usepackage{amssymb}

\usepackage[text={130mm,200mm}]{geometry}

\theoremstyle{plain}

\newtheorem*{theorem*}{Theorem}

\newtheorem*{corollary*}{Corollary}

\newtheorem*{lemma*}{Lemma}

\newtheorem*{proposition*}{Proposition}

\newtheorem*{conjecture*}{Conjecture}
\theoremstyle{definition}

\newtheorem*{definition*}{Definition}
\theoremstyle{remark}

\newtheorem*{remark*}{Remark}
\newtheorem{example}{Example}

\providecommand{\surname}[1]{#1}
\providecommand{\country}[1]{#1}

\begin{document}

\title[Finite shtukas]{On  Langlands program, global fields and shtukas}

\author{Nikolaj \surname{Glazunov}}
\address{Department of Electronics\\
  National Aviation University\\
  1~Komarova Pr.\\
  Kiev\\
 03680 \\
  \country{Ukraine}}
\email{glanm@yahoo.com}
\urladdr{https://sites.google.com/site/glazunovnm/}



\date{June 25, 2020}

\begin{abstract}
The purpose of this paper is to survey some of the important results on 
Langlands program, global fields, $D$-shtukas and
finite shtukas which have influenced the development of  algebra and number theory. It is intended to be 
selective rather than exhaustive, as befits the occasion of the 80-th birthday of Yakovlev, 75-th birthday of Vostokov and 
75-th birthday of Lurie. \\
Under assumptions on ground fields results on Langlands program have been proved and discussed by
Langlands, Jacquet, Shafarevich, Parshin, Drinfeld, Lafforgue and others.

This communication is an introduction to the Langlands Program, global fields
and to  $D$-shtukas and
finite shtukas  (over algebraic curves) over function fields.
At first recall that linear algebraic groups found important applications in the Langlands program. Namely, for a connected reductive group $ G $ 
over a global field $ K $, the Langlands correspondence relates automorphic forms on $ G $ and global 
Langlands parameters,  i.e. conjugacy classes of homomorphisms  from the Galois group
 $ {\mathcal Gal} ({\overline K} / K) $ 
to the dual Langlands group $ \hat G ({\overline {\mathbb Q}} _ p) $. In the case of fields of 
algebraic numbers, the application and development of elements of the Langlands program made it 
possible to strengthen  the Wiles theorem on the Shimura-Taniyama-Weil hypothesis and to prove
 the Sato-Tate hypothesis.
 
  V. Drinfeld and L. Lafforgue have investigated the case of functional global fields of characteristic $ p> 0 $
   ( V. Drinfeld for $ G = GL_2 $ and L. Lafforgue for $ G = GL_r, \; r $ is an arbitrary positive integer). 
   They have proved in these cases the Langlands correspondence.
 
 Under the process of these investigations, V. Drinfeld  introduced the concept of a $ F $-bundle, or shtuka,
 which was used by both authors in the proof for functional global fields of characteristic $ p> 0 $ 
of the studied cases of the existence of the Langlands correspondence.
 
 Along with the use of shtukas developed and used by V. Drinfeld and L. Lafforge, other constructions
 related to approaches to the Langlands program in the functional case were introduced.
 
 G. Anderson has introduced the concept of a $ t $-motive. U. Hartl, his colleagues and students have
 introduced and have explored the concepts of finite, local and global $ G $-shtukas.
 
 In this review article, we first present results on Langlands program and  related representation over algebraic number fields.
 Then we  briefly present approaches by U. Hartl, his colleagues and students to the study
 of $ D $ --shtukas and finite shtukas.
 These approaches and our discussion relate to the Langlands program as well as to the internal development
  of the theory of $ G $-shtukas.

\end{abstract}

\maketitle


\section*{Introduction}

This communication is an introduction to the Langlands Program
and to ($D$-)shtukas and finite shtukas (over algebraic curves) over function fields.
 The Langlands correspondence over number fields in its full generality is facing with 
problems \cite{Langlands1979,Lan1980,JacLan1970,Sha96,Parshin,Drinfeld1980,Lafforgue2002}. So results 
 from Galois theory,  algebraic number theory and function fields can help understand it.
\subsection{Elements of algebraic number theory and field theory.}
The questions what is a Galois group of a given algebraic closure of  the number field or the local field, embedding problems of 
fields and extensions of class field theory belong to fundumental questions of Galois theory and class field theory.
A.V. Yakovlev,  S.V. Vostokov,  B.B. Lur'e works spans many areas of Galois theory, fields theory
and class field theory.
The results obtained
indicate that these questions  connect with module theory, homological algebra and with other topics of algebra and number theory\cite{Sha,Yak1967,Yak,Vos11,Vos85,Lur1964,Lur1991}.
The development and applications of these theories are discribed  in papers by I.R. Shafarevich\cite{Sha96} and by F.N. Parshin \cite{Parshin} 
(and in references therein). 
For further details we refer the reader to papers themselves.
By the lack of author`s competence we discuss here very shortly only  connection of local fields with formal modules.

\subsection{The Hensel-Shafarevich canonical basis in complete discrete valuation fields.}
Vostokov has constructed a canonical Hensel-Shafarevich basis in  ${\mathbb Z}_p-$module of principle units for complete discrete valuation field with an arbitrary residue field \cite{Vos11}.
Vostokov and Klimovitski in paper \cite{VosKl}  give construction of primary elements in formal modulei.
Ikonnikova, Shaverdova \cite{IKSh} and Ikonnikova \cite{IK} use these results  under  construction, respectively,
 the Shafarevich basis in higher-dimensional
 local fields   and under proving two theorems  on the canonical 
basis in Lubin-Tate formal modules in 
the case of local field with perfect residue field  and in the case of imperfect 
residue field. 
These canonical bases are obtained by applying a variant of the Artin-Hasse function. 

\subsection{${^L} G$ for  reductive group $G$}
Here we follow to  \cite{Langlands1979,JacLan1970,Arthur,Borel,Tate}.
At first recall that linear algebraic groups found important applications in the Langlands program.
 Namely, for a connected reductive group $ G $ 
over a global field $ K $, the Langlands correspondence relates automorphic forms on $ G $ and global 
Langlands parameters,  i.e. conjugacy classes of homomorphisms  from the Galois group
 $ {\mathcal Gal} ({\overline K} / K) $ 
to the dual Langlands group $ \hat G ({\overline {\mathbb Q}} _ p) $.
Let $\overline K$ be an algebraic closure of $ K$ and $ K_s$  be the separable  closure of $ K $ in $\overline K$.
\begin{definition*}
Let $G$ be the connected reductive algebraic group over $\overline K$. 
The root datum of $G$ is a quadruple $(X^*(T), \Delta, X_*(T), \Delta^v)$   where $X^*$ is the lattice of characters of the maximal torus $T$, 
$ X_*$ is the dual lattice, given by the 1-parameter subgroups, $\Delta$ is the set of roots, $\Delta^v$  is the corresponding set of coroots.
\end{definition*}

The dual Langlands group $ \hat G$ is a complex reductive group that has the dual root data:  $(X_*(T), \Delta^v, X^*(T), \Delta )$.
Here any maximal torus  $ \hat T$ of  $ \hat G$ is isomorphic to the complex dual torus 
$ X^*(T) \otimes {\mathbb C}^* = Hom(X_*(T),{\mathbb C}^*) $ of any maximal torus $T$ in $G$. 
Let $\Gamma_{\mathbb Q} = {\mathcal Gal} ({\overline Q} / Q) $.

Given $G$, the Langlands ${L}$-group of $ G$ is defined as semidirect product
\begin{equation*}
{^L} G =  \hat G \rtimes \Gamma_{\mathbb Q}.
\end{equation*}
 
 In the case of fields of 
algebraic numbers, the application and development of elements of the Langlands program made it 
possible to strengthen
 the Wiles theorem on the Shimura-Taniyama-Weil hypothesis and to prove the Sato-Tate hypothesis.
 Langlands reciprocity for $GL_n$ over non-archimedean local fields of characteristic zero is given by 
 Harris-Taylor \cite{HarrisTaylor}.

\subsection{ Langlands correspondence over functional global fields of characteristic $ p> 0 $}
 V. Drinfeld \cite{Drinfeld1980} and L. Lafforgue \cite{Lafforgue2002} have investigated the case of functional global fields of characteristic $ p> 0 $
   ( V. Drinfeld for $ G = GL_2 $ and L. Lafforgue for $ G = GL_r, \; r $ is an arbitrary positive integer). 
   They have proved in these cases the Langlands correspondence.
 
 In the process of these studies, V. Drinfeld introduced the concept of a $ F $ -bundle, or shtuka,
 which was used by both authors in the proof for functional global fields of characteristic $ p> 0 $ 
of the studied cases of the existence of the Langlands correspondence \cite{Drinfeld1987}.
 
 Along with the use of shtukas developed and used by V. Drinfeld and L. Lafforge, other constructions
 related to approaches to the Langlands program in the functional case were introduced.
 
 G. Anderson has introduced the concept of a $ t $-motive \cite{Anderson1986}. U. Hartl, his colleagues and postdoc students have
 introduced and have explored the concepts of finite, local and global $ G $ -shtukas\cite{Hartl2005,HartlViehmann2011,HartlRad2014,Singh2012,Rad2012,Weiss2017}.
 
 In this review, we first present results on Langlands program and  related representation over algebraic number fields.
 Then we  briefly present approaches by U. Hartl, his colleagues and students to the study
 of $ G $ --shtukas.
 These approaches and our discussion relate to the Langlands program as well as to the internal development
  of the theory of $ G $-shtukas.
  Some results on commutative formal groups and commutative formal schemes can be found in  
   \cite{Glazunov2015A,Glazunov2015B,Glazunov2015C} and in references therein.
  
 The content of the paper is as follows: \\
   Introduction. \\
  1. Some results of the implementation of the Langlands program for fields of algebraic numbers and their localizations. \\
  2. Elliptic modules and Drinfeld shtukas. \\
  3. Finite $ G $ -shtukas. \\

\section{Some results on Langlands program over algebraic number fields and their localizations} 

    Langlands conjectured that some symmetric power $ L$-functions extend to an entire function and coincide with certain automorphic $L$-functions.
    
\subsection{ Abelian extensions of number fields }

     In the case of algebraic number fields Langlands conjecture (Langlands correspondence) is the  global class field theory:    \\
          Representations of the abelian Galois group $Gal(K^{ab}/K)$ = characters of the Galois group $Gal(K^{ab}/K) $    
\begin{center}
                                                     correspond to  
\end{center}
automorphic forms  on $GL_1 $ that are characters of the class group of ideles. Galois group    $ Gal(K^{ab}/K)$ is the profinite completion of the group ${\mathbb A}^* (K)/K^* $ where $ {\mathbb A}(K) $ denotes the adele ring of $ K$.\\
  If $K$ is the local field, then Galois group $ Gal(K^{ab}/K)$ is canonically isomorphic to the profinite completion of $K^*$.   
    
\subsection{ $l$ - adic representations and Tate modules}

       Let $K$ be a field and $\overline K $  its separate closure, $ E_n=\{P \in E(\overline K )| nP=0\}$ the group of points of elliptic curve $E(\overline K)$  order dividing n. When $char K $ does not divide $n$ then $E_n$ is a free ${\mathbb Z}/n{\mathbb Z}$ -module of rank $2$.
       
             Let $l$ be prime, $l \ne char K$. The projective limit $T_l (E)$ of the projective system of modules $E_{l^m}$ is free
${\mathbb Z}_l$-adic Tate module of rank $2$.  

Let              
$V_l (E)=T_l (E) \otimes_{{\mathbb Z}_l} {\mathbb Q}_l$.
      Galois group $Gal({\overline K}/K)$ acts on all $E_{l^m}$ , so there is the natural continuous representation       ($l$-adic representation)
\begin{equation*}
\rho_{E,l} : Gal({\overline K}/K) \to Aut \; T_l (E) \subseteq Aut \;  V_l (E).  
\end{equation*}
$V_l (E)$ is the first homology group that is dual to the first cohomology group of $ l$-adic cohomology of elliptic curve $E$ and Frobenius $F$ acts on the homology and dually on cohomology. The characteristic polynomial $ P(T)$ of the Frobenius not depends on the prime number $ l$.

\subsection{Zeta functions and parabolic forms}

Let (in P. Deligne notations) $X$ be a scheme of finite type over ${\mathbb Z}$, $|X|$ the set of its closed points, and for each $x \in |X|$ let  $N(x)$ be the number of points of the residue field $k(x)$ of $X$ at $x$. The Hasse-Weil zeta-function of $X$ is, by definition
\begin{equation*}
\zeta_X (s)= \prod_{ x \in|X|}(1 - N(x)^{-s} )^{-1} .
\end{equation*}
In the case when   $X$  is defined over finite field ${\mathbb F}_q$, put  $ q_x= N(x)$, $deg(x)=[k(x):{\mathbb F}_q ]$, so $q_x=q^{deg(x)}$ . Put $t=q^{-s}$. Then 
\begin{equation*}
Z(X ,t) = \prod_{ x \in|X|}(1 - t^{deg(x)} )^{-1} .
\end{equation*}

The Hasse-Weil zeta function of $ E$ over ${\mathbb Q}$ (an extension of numerators of $\zeta_E (s)$ by points of bad reduction of $E$) is defined over all primes $p$:
\begin{equation*}
L(E({\mathbb Q}),s) = \prod_p (1 - a_pp^{-s} + \epsilon(p)p^{1 - 2s} )^{-1} ,
\end{equation*}
here $\epsilon(p) = 1$ if $E$ has good reduction at $ p$, and $\epsilon(p) = 0$ otherwise.

Put $T = p^{-s}$. For points of good reduction we have 
\begin{equation*}
P(T) = 1 - a_pT + pT^2  = (1- \alpha T)(1 - \beta T)
\end{equation*}.
For symmetric power $L$-functions (functions $  L(s;E; Sym^n), \; n > 0;$ see below) we have to put
\begin{equation*}
P_p(T) = \prod_{i=0}^n (1- \alpha^i \beta^{n - i}T)
\end{equation*}

For $GL_2 ({\mathbb R})$, let $C$ be its center,  $O(2)$ the orthogonal group.

Upper half complex plane has the representation: ${\mathbb H}^2=GL_2 ({\mathbb R}) / O(2) C$. So it is the homogeneous space of the group  $GL_2 ({\mathbb R})$.

A cusp (parabolic) form of weight $k \ge 1$ and level $N \ge1$ is a holomorphic function $f$ on the upper half complex plane ${\mathbb H}^2$ such that \\
a) For all matrices

$$g = \left(
\begin{array}{cc}
a& b \\
c& d
\end{array}
\right),     a,b,c,d \in {\mathbb Z}, a \equiv1(N),d\equiv1(N),c\equiv0(N) $$

and for all  $z \in {\mathbb H}^2$ we have
\begin{equation*}
f(gz)=f((az+b)/(cz+d))=(cz+d)^k f(z)
\end{equation*}
(automorphic condition).\\
b) 
\begin{equation*}
	|f(z)|^2 (Im z)^k 
\end{equation*}
is bounded on ${\mathbb H}^2$ . \\

 Mellin transform $ L(f,s)$ of the parabolic form $ f$ coincides with Artin $ L$-series of the representation $\rho_f$.
 
 The space ${\mathcal M}_n(N)$  of cusp forms of weight $ k$ and level  $N$ is a finite dimensional complex vector space. If 
 $f \in {\mathcal M}_n (N)$, then it has expansion 
 \begin{equation*}
f(z)=\sum_{n=1}^{\infty} c_n (f)\exp(2\pi inz)
\end{equation*}
and $L$-function is defined by
 \begin{equation*}
 L(f,s)=\sum_{n=1}^{\infty} c_n (f) /  n^s .
\end{equation*}

\subsection{Modularity results} 
  
  The compact Riemann surface $\Gamma\backslash{\mathbb H}^2$ is called the modular curve associated to the subgroup of finite index 
$\Gamma$ of  $GL_2 ({\mathbb Z})$ and is denoted by $X(\Gamma)$. If the modular curve is elliptic it is called the elliptic modular curve. \\    
The modularity theorem states that any elliptic curve over ${\mathbb Q}$ can be obtained via a rational map with integer coefficients from the elliptic  modular curve. \\
      By the Hasse-Weil conjecture (a cusp form of weight two and level $ N$ is an eigenform (an eigenfunction of all Hecke operators)). The conjecture follows from the modularity theorem.\\
         Recall the main (and more stronger than in Wiles \cite{Wiles} and in Wiles-Taylor \cite{TaylorWiles} papers) result by C. Breuil, B. Conrad, F. Diamond, R. Taylor \cite{BreuilConradDiamondTaylor}.
 
\begin{theorem*}
({\bf Taniyama-Shimura-Weil conjecture - Wiles Theorem.})
For every elliptic curve $E$ over ${\mathbb Q}$ there exists $f$, a cusp form of weight 2 for a subgroup 
$\Gamma_0 (N)$, 
such that $ L(f,s)=L(E({\mathbb Q}),s)$. 
\end{theorem*}
Here ${\Gamma_0(N)}$ is the modular group 
\begin{equation*}
{\Gamma_0(N)} = \left\{
\left(
  \begin{array}{cc} 
  a& b\\
  c& d \\
   \end{array}
   \right),
     a,b,c,d \in {\mathbb Z}, c \equiv 0\pmod{N},
     \det \left(
  \begin{array}{cc} 
  a& b\\
  c& d \\
   \end{array}
   \right)= 1 \right \}.
\end{equation*}
Recall that for projective closure $\overline E$ of the elliptic curve $E$ we have
\begin{equation*}
\overline E({\mathbb F}_p) = 1 - a_p +p. 
\end{equation*}
By H. Hasse
\begin{equation*}
a_p = 2\sqrt p \cos \varphi_p.
\end{equation*}
\begin{conjecture*}
({\bf Sato-Tate conjecture})
           Let $E$ be an  elliptic curve  without complex multiplication.
             Sato have computed and Tate gave theoretical evidence that angles $\varphi_p$ in the case are equidistributed in $[0,\pi]$ with the Sato-Tate density measure  $\frac{2}{\pi} \sin^2 \varphi.$
\end{conjecture*}
 
We have  two theorems from  Serre  \cite{Serre} which give the theoretical explanation in terms of Galois representations.
Here we recall the corollery of the theorems.
\begin{corollary*}
(Serre  \cite{Serre})
The elements 
are equidistributed for the $v$ normalized Haar measure of $G$ if and 
only if $c = 0$ for every $X$ irreducible character 
of $G$, i. e. , if and only if the 
$L$-functions relative to the non trivial irreducible characters of $G$ are holomorphic and 
non zero at $s = 1$. 
\end{corollary*}

The current state of Sato-Tate conjecture is now Clozel--Harris--Shepherd-Barron--Taylor Theorem 
\cite{Clozel HarrisTaylor,HarrisShepherd-BarronTaylor}.

\begin{theorem*}
(Clozel, Harris, Shepherd-Barron, Taylor). Suppose $E$ is an elliptic curve over ${\mathbb Q}$ with non-integral $j$ invariant. Then for all $n > 0; \; L(s;E; Sym^n)$ extends to a meromorphic function which is holomorphic and non-vanishing for $Re(s) \ge 1 + n/2$.
\end{theorem*}

These conditions and statements are sufficient to prove the Sato-Tate conjecture.

Under the prove of the Sato-Tate conjecture the Taniyama-Shimura-Weil conjecture oriented methods of A. Wiles and R. Taylor are used.

     Recall also that the proof of Langlands reciprocity for $GL_n$ over non-archimedean local fields of characteristic zero is given by Harris-Taylor \cite{HarrisTaylor}.
     
\section{Elliptic modules and Drinfeld shtukas.}

Let \\
${\overline{\mathbb F}}_q$ be the algebraic closure of ${\mathbb F}_q$, \\
$\mathcal C$  be  a smooth projective geometrically irreducible curve over ${\mathbb F}_q$, \\
$K$ be  the function field  ${\mathbb F}_q (\mathcal C)$ of $\mathcal C$, \\
 $\nu$ be a close point of  $\mathcal C$,  \\
$A$ be the ring of functions regular on  $\mathcal C - \nu$, \\
$K_{\nu}$ be the complation of $K$ at $\nu$ with valuation ring ${\mathcal O}_{\nu}$, \\
${\mathbb C}_{\nu}$ be the complation of  the algebraic closure of $K_{\nu}$.  \\
At first recall some known facts about algebraic curves over finite fields.
We will identify   the set $|\mathcal C|$ of closed points of  $\mathcal C$ with 
$\mathcal C({\overline{\mathbb F}}_q) = Hom_{{\mathbb F}_q} (Spec \; {\overline{\mathbb F}}_q, \mathcal C)$.
Let $k(\nu)$ be the residue field of $\nu$. Then the degree of $\nu$ is equal of the number of elements $[k(\nu):{\mathbb F}_q]$.

Below in this section we follow to \cite{Drinfeld1974,Drinfeld1987,DeHu}.
\subsection{Elliptic modules}
\begin{lemma*}
Let $k$ be a field of characteristic $p > 0$ and let $R$ be a $k$-commutative ring with unit (there exists a morphism $k \to R$).  The additive scheme ${\mathbb G}_a$ over $R$ is represented by the polynomial ring $R[X] $ with structural morphism 
$\alpha:R[X] \to R[X]\otimes_R R[X] $, given by $ \alpha(X) = X \otimes 1 + 1 \otimes X$.
A  morphism $\varphi: {\mathbb G}_a \to {\mathbb G}_a$ of additive schemes over $R$ is defined by an additive polynomial.
If $\psi$ is another such morphism, then $\varphi\circ\psi = \varphi(\psi(T))$. So the set of (endo)morphisms of additive scheme has the structure of a ring.
\end{lemma*}
\begin{example}
Let $a \in R[X], \; pa = 0$. Then the morphism $\varphi(T) = a T^{p^n}, \; (n \ge 0)$ is additive.
Any additive morphism  $\varphi(T)$ in characteristic $p$ has the form 
$\varphi(T) = a_0T + a_1T^p + \cdots + a_n T^{p^n}.$
\end{example}
\begin{proposition*}
Let $k$ be a field of characteristic $p > 0$.
Put $\tau a = a^p\tau$. There is an isomorphism between $End_k({\mathbb G}_a)$ and the ring  of noncommutative 
polynomials $k\{\tau\}$.
\end{proposition*}
For any $\varphi(T) = a_0T + a_1T^p + \cdots + a_n T^{p^n} \in End_k({\mathbb G}_a)$ and any 
$ \varphi(\tau) = a_0 + a_1\tau + \cdots + a_n\tau^n \in k\{\tau\}$ Lubin morphisms \cite{Lub} $c_0$ and $c$ are defined:
\begin{equation*}
c_0(\varphi(T)) = a_0,   c(\varphi(\tau)) = a_0.
\end{equation*}
Respectively we  define 
\begin{equation*}
deg(\varphi(T)) = p^n,  d(\varphi(\tau)) = n.
\end{equation*}
\begin{proposition*}
Any ring morphism $A \to End_k({\mathbb G}_a)$  is either injective or has image contained in the constants $k \subset k\{\tau\}$.
\end{proposition*}
Sketch of the proof. $k\{\tau\}$ is a domain. $End_k({\mathbb G}_a)$ is isomorphic to $k\{\tau\}$.
$A$ is a ring with divisor theory ${\mathfrak D}$ and for any prime divisor 
${\mathfrak p} \in {\mathfrak D}$ the residue ring
$A/{\mathfrak p}$ is a field. From these statements the proposition follows. \\
Assume now that $k$ is an $A$-algebra, i.e. there is a morphism $i: A \to k$ .
\begin{definition*}
An elliptic module over $k$ (of rank $r =2$) is an injective ring homomorphism 
\begin{equation*}
 \varphi: A \to End_k({\mathbb G}_a)
\end{equation*}
\begin{equation*}
  a \mapsto \varphi_a ,
\end{equation*}
 such that for all $a \in A$ we have 
\begin{equation*}
  d(\varphi(\tau)) = 2\cdot deg(a),
 \end{equation*}
 \begin{equation*}
   c(\varphi(\tau)) = i(a).
\end{equation*}
\end{definition*}

\begin{example}
 Let $k = {\mathbb F}_q(T)$, $A = {\mathbb F}_q[{\mathbb P}^1 - \nu] = {\mathbb F}_q[T]$.
 Let $i(T) = T^2 +1$. 
 In this case an elliptic module $\varphi$ is given by 
 \begin{equation*}
   \varphi = T^2 +1 + c_1\cdot\tau + c_2\cdot\tau^2, c_1, c_2 \in k, \; c_2 \ne 0.
\end{equation*}
\end{example}
\begin{remark*}
By the same way it is possible to define a Drinfeld module (over a field) for any natural $r$.
\end{remark*}
Now consider the case of Drinfeld modules over a base scheme. Let $S$ be an $A$-scheme,
 $\mathcal L$ a line bundle over $S$, $i^*: S \to Spec \; A$ be an $A$ scheme morphism dual to the ring homomorphism $i:A \to  {\mathcal O}_S$
\begin{definition*}
(Drinfeld module over a base scheme)
A Drinfeld module over $k$ of rank $r$ is an ring homomorphism 
\begin{equation*}
 \varphi: A \to End_S(\mathcal L)
\end{equation*}
\begin{equation*}
  a \mapsto \varphi_a ,
\end{equation*}
 such that for all $a \in A$ we have \\ 
1)  locally, as a polynomial in $\tau$, $\varphi_a$ has the degree 
\begin{equation*}
  d(\varphi(\tau)) = r\cdot deg(a),
 \end{equation*}
 2) a unit as its leading coefficient $a_n$ and
  \begin{equation*}
  c(\varphi(\tau)) = i(a).
\end{equation*}
\end{definition*}
\subsection{Drinfeld shtukas.}
In notations of previous subsection let $x \in k$, $a \in A$, $\varphi_a(\tau)$ be a Drinfeld module of rank $r$. 
Put $L =  k\{\tau\}$, $f(\tau) \in L$, $k[A] = k\otimes_{{\mathbb F}_q}A$, $deg_{\tau}f(\tau)$ the degree in $\tau$ of $f(\tau)$.
\begin{lemma*}
Define the action  of $k[A]$ on $L$ by the formula:
 \begin{equation*}
  x\otimes a\cdot f(\tau) = x\cdot f(\varphi_a(\tau) .
\end{equation*}
Then $L$  is a free $k[A]$-module of rank $r$.
\end{lemma*}
\begin{remark*}
Let $E_s = \{f(\tau)  \in L | deg_{\tau}f(\tau) \le s\}, \;  E = \oplus_{s=0}^\infty E_s$,
$ E[1] = \oplus_{s=0}^\infty E_{s+1}.$
$ E,  E[1]$ are graded modules over the graded ring and give rise to locally free sheaves ${\mathcal F}$, ${\mathcal E}$ of rank $r$ over ${\mathcal C}$.
\end{remark*}
Put ${\mathcal C}_S = C \times_{{\mathbb F}_q}S$, 
$\sigma_q = id_C \otimes Frob_{q,S}:{\mathcal C}_S \to {\mathcal C}_S $
\begin{definition*}
A (right) $\mathcal D$-shtuka ($F$-sheaf \cite{Drinfeld1987})  of rank $r$ over an 
${{\mathbb F}_q}$-scheme $S$ is a diagram 
$({\mathcal F} \stackrel{c_1}\to {\mathcal E} \stackrel{c_2}\gets ( id_C \otimes Frob_{q,S})^*{\mathcal F}  )$, such that
coker $c_1$ is supported on the graph ${\Gamma}_{\alpha}$ of a morphism ${\alpha}: S \to {\mathcal C}$
 and it is a line bundle on support,
 coker $c_2$ is supported on the graph ${\Gamma}_{\beta}$ of  a morphism ${\beta}: S \to {\mathcal C}$
 and it is a line bundle on support.
 \end{definition*}
If  ${\Gamma}_{\alpha} \cap {\Gamma}_{\beta} = \emptyset$ it is possible to give the next
definition of $\mathcal D$-shtuka \cite{HartlRad2014,HartlSingh}.
\begin{definition*}
A  global shtuka  of rank $r$ with two legs over an 
${{\mathbb F}_q}$-scheme $S$ is a tuple 
$\underline{\mathcal N} =  ({\mathcal N}, (c_1, c_2), {\tau}_{\mathcal N})$ 
consisting of 1) a locally free sheaf ${\mathcal N}$ of rank $r$ on ${\mathcal C}_S$;
2) ${{\mathbb F}_q}$-morphisms $c_i: S \to {\mathcal C} \; ( i=1,2)$, called the legs of 
$\underline{\mathcal N}$; 3) an isomorphism 
$\tau_N: \sigma_q^*{\mathcal N}|_{{\mathcal C}_S - ({\Gamma}_{c_1} \cup {\Gamma}_{c_2})}\simeq {\mathcal N}|_{{\mathcal C}_S - ({\Gamma}_{c_1} \cup {\Gamma}_{c_2})}$
outside the graphs ${\Gamma}_{c_i}$  of $c_i$, ${\Gamma}_{c_1} \cap {\Gamma}_{c_2} = \emptyset$.
\end{definition*}
\begin{definition*}
A  global shtuka over $S$ is a  $\mathcal D$-shtuka if $\tau_N$ satisfies
$
\tau_N(\sigma_q^*{\mathcal N}) \subset {\mathcal N}
$
on ${\mathcal C}_S - {\Gamma}_{c_2}$ with cokernel locally free of rank 1 as $\mathcal O_S$-module, and 
$\tau_{\mathcal N}^{-1}({\mathcal N}) \subset \sigma_q^*{\mathcal N}$ on ${\mathcal C}_S - {\Gamma}_{c_1}$ with cokernel locally free of rank 1 as $\mathcal O_S$-module.
\end{definition*}

\section{Finite $ G $-shtukas.}
We follow to \cite{Drinfeld1987,Hartl2005,HartlViehmann2011,Singh2012,HartlSingh}.
We start with very short indication on the general framework of the section. In connection with Drinfeld`s constructions of elliptic  modules Anderson \cite{Anderson1986}  has introduced abelian t-modules and the dual notion of t-motives. Beside with mentioned papers these are the descent theory by A. Grothendieck \cite{Grothendieck}, cotangent complexes  by Illusie \cite{Illusie}, by S. Lichtenbaum and M. Schlessinger \cite{LichtenbaumSchlessinger}, by Messing \cite{Messing}  and by Abrashkin \cite{Abrashkin}. In this framework to any morphism $f: A \to B$ of commutative ring objects in a topos is associated a cotangent complex $ L_{(B/A)}$ and to any morphism of commutative ring objects in a topos of finite and locally free 
$Spec (A)$-group schemes $G$ is associated a cotangent complex $L_{(G/Spec (A))}$  as has presented in books by  Illusie \cite{Illusie}.
\subsection{Finite shtukas and formal groups} 
Let $S$ be a scheme over $Spec \; {\mathbb F}_q$. 
\begin{definition*}
A finite ${\mathbb F}_{q}$-shtuka over $S$ is a pair ${\underline M} = (M, F_M )$ consisting of a locally
 free  ${\mathcal O}_S$-module $M$ on $S$ of finite rank  and an ${\mathcal O}_S$-module homomorphism 
 $F_M : \sigma_q^*M \to M$.
\end{definition*}
 Author\cite{Singh2012} investigates relation between finite shtukas and strict finite flat commutative group schemes and relation between divisible local Anderson modules and formal Lie groups.  
 The cotangent complexes as in papers by S. Lichtenbaum and M. Schlessinger  
\cite{LichtenbaumSchlessinger}, by W. Messing \cite{Messing}, by V. Abrashkin 
\cite{Abrashkin} are defined and are proved that they are homotopically equivalent. \\
  Then the deformations of affine group schemes follow to the mentioned paper of Abrashkin are investigated   and strict finite $\mathcal O$-module schemes are  defined.\\
Next step of the research  is devoted to relation between finite shtukas by V. Drinfeld \cite{Drinfeld1987} and strict finite flat commutative group schemes.\\
   The comparison between cotangent complex and Frobenius map of finite ${\mathbb F}_p$-shtukas is given.\\
   \subsection{Local shtukas and local Anderson modules} 
   Recall some notions and notations. An ideal $I$ in a commutative ring $A$ is locally nilpotent at a prime ideal $\varrho$ if the 
localization $I_{\varrho}$ is a nilpotent ideal in $A_{\varrho}$.
In the framework of   smooth projective geometrically irreducible curves $\mathcal C$ over ${\mathbb F}_q$ let 
$Nilp_{{A}_{\nu}}$ denote the category of  $A_{\nu}$-schemes on which the uniformizer $\xi$ of $A_{\nu}$ is locally
 nilpotent. Here $A_{\nu} \simeq {\mathbb F}_{\nu}[[\xi]]$ is the completion of the local ring $\mathcal O_{\mathcal C, \nu}$
 at a closed point $\nu \in \mathcal C$. \\
 Let $Nilp_{{\mathbb F}_{q}[[\xi]]}$ be the category of ${\mathbb F}_{q}[[\xi]]$-schemes on which $\xi$ is locally nilpotent. Let $S \in Nilp_{{\mathbb F}_{q}[[\xi]]}$. Let $M$ be  a sheaf of $\mathcal O_S[[z]]$-modules on $S$ and let
$\sigma_q^*M = M \otimes_{\mathcal O_S[[z]],\sigma_q^*} O_S[[z]]$,
 $M[\frac{1}{z - \xi}] = M\otimes_{\mathcal O_S[[z]]} {\mathcal O_S[[z]]}[\frac{1}{z - \xi}]$ . 
 \begin{definition*}
  A local shtuka of  height $r$ over $S$ is a pair $M = (M, F_M )$ consisting of a
locally free sheaf $M$ of ${\mathcal O}_S[[z]]$-modules of rank $r$, and an isomorphism
 $F_M: \sigma_q^*M[\frac{1}{z - \xi}] \simeq M[\frac{1}{z - \xi}]$.
 \end{definition*}
 The next lemma is proved \cite{HartlSingh}.
 \begin{lemma*}
  Let $R$ be an ${\mathbb F}_{q}[[\xi]]$-algebra in which $\xi$ is nilpotent. Then the sequence of $R[[z]]$-modules
   \begin{equation*}
   0 \to R[[z]] \to R[[z]] \to R \to 0
    \end{equation*}
    \begin{equation*}
      1 \mapsto z- \xi,   z \mapsto  \xi
    \end{equation*}
    is exact. In particular $R[[z]] \subset R[[z]][\frac{1}{z - \xi}]$.
 \end{lemma*}
 In the conditions of the lemma authors \cite{HartlSingh} give the next 
  \begin{definition*}
 A $z$-divisible local Anderson module over $R$ is a sheaf of ${\mathbb F}_{q}[[z]]$-modules $G$ on the big
$fppf$-site of $Spec \; R$ such that\\
(a) $G$ is $z$-torsion, that is $G = \varinjlim  G[z^n]$, where $G[z^n] = ker(z^n: G \to G)$, \\
(b) $G$ is $z$-divisible, that is $z : G \to G$ is an epimorphism, \\
(c) For every $n$ the ${\mathbb F}_{q}$-module $G[z^n]$  is representable by a finite locally free strict 
${\mathbb F}_{q}$-module scheme
over $R$ in the sense of Faltings (\cite{Faltings, HartlSingh}), and \\
 (d) locally on $Spec \; R$ there exists an integer $ d \in {\mathbb Z}_{\ge 0}$, such that 
 $ (z - \xi)^d = 0$ on $\omega_G$ where $\omega_G  = \varprojlim \omega_{G[z^n]}$ and 
 $\omega_{G[z^n]} = \varepsilon^* {\Omega^1}_{G[z^n]/{Spec \; R}}$ for the unit section 
  $\varepsilon$ of $G[z^n]$ over $R$.
   \end{definition*}
    $ z$-divisible local Anderson modules by  Hartl \cite{Hartl2005} with improvements in \cite{HartlSingh}  and local shtukas are investigated.
   The equivalence between the category of effective local shtukas over $S$ and the category of $z$-divisible local Anderson modules over $S$ is treated by the authors \cite{Singh2012,HartlSingh}. 
   The theorem about canonical ${\mathbb F}_p[[\xi]]$ -isomorphism of $z$-adic Tate-module of $z$ -divisible local Anderson module $G$ of rank $r$ over $S$ and Tate module of local shtuka over $S$ associated to $G$ is given.
   The main result of \cite{Singh2012} is the following (section 2.5) interesting result:
 it is possible to associate a formal Lie group to any $z$-divisible local Anderson module over $S$ in the case when $\xi$ is locally nilpotent on $S$. We note that related with \cite{Singh2012} and in some cases more general results have 
presented in the paper by U. Hartl, E. Viehmann \cite{HartlViehmann2011}.

\newpage
\end{document}